\RequirePackage{ifpdf}
\ifpdf 
\documentclass[pdftex]{sigma}
\else
\documentclass{sigma}
\fi

\newcommand{\nc}{\newcommand}

\nc{\slt}{\mathfrak{sl}_2}
\nc{\suth}{\widehat{\mathfrak{su}}(2)}
\nc{\gl}{\mathfrak{gl}}
\nc{\GL}{\mathfrak{GL}}
\nc{\g}{\mathfrak{g}}
\nc{\alg}{\mathfrak{a}}
\nc{\gh}{\widehat\g}
\nc{\h}{\mathfrak{h}}
\nc{\la}{\lambda}
\nc{\slth}{\widehat{\slt}}
\nc{\C}{\mathbb C }
\nc{\Z}{\mathbb Z }
\nc{\N}{\mathbb N }
\nc{\R}{\mathbb R }
\nc{\Q}{\mathbb Q }
\nc{\al}{\alpha }
\nc{\be}{\beta}
\nc{\ta}{\theta}
\nc{\ve}{\varepsilon}
\nc{\ch}{{\mathop {\rm ch}}}
\nc{\Tr}{{\mathop {\rm Tr}\,}}
\nc{\Id}{{\mathop {\rm Id}}}
\nc{\U}{{\mathop {\rm U}}}
\nc{\bra}{\langle}
\nc{\ket}{\rangle}
\nc{\x}{{\bf x}}
\nc{\pa}{\partial}
\nc{\ld}{\ldots}
\nc{\cd}{\cdots}
\nc{\hk}{\hookrightarrow}
\nc{\A}{\mathfrak A}
\nc{\qb}[2]{\genfrac{(}{)}{0pt}{}{#1}{#2}_q}
\nc{\n}{\mathfrak{n}}
\nc{\un}{\mathfrak{u}}
\nc{\T}{\otimes}
\nc{\bv}{{\bf v}}
\nc{\bu}{{\bf u}}
\nc{\bs}{{\bf s}}
\nc{\bw}{{\bf w}}
\nc{\qbin}[2]{{\genfrac{(}{)}{0pt}{0}{#1}{#2}}_q}
\nc{\fac}[1]{(#1)_q!}
\nc{\wt}{\widetilde}
\nc{\at}[2]{\genfrac{}{}{0pt}{}{#1}{#2}}
\nc{\gr}{\mathrm{gr}}
\nc{\ov}{\overline}

\begin{document}


\renewcommand{\thefootnote}{$\star$}

\renewcommand{\PaperNumber}{070}

\FirstPageHeading

\ShortArticleName{The PBW Filtration, Demazure Modules and Toroidal Current Algebras}

\ArticleName{The PBW Filtration, Demazure Modules\\ and Toroidal Current Algebras\footnote{This paper is a
contribution to the Special Issue on Kac--Moody Algebras and Applications. The
full collection is available at
\href{http://www.emis.de/journals/SIGMA/Kac-Moody_algebras.html}{http://www.emis.de/journals/SIGMA/Kac-Moody{\_}algebras.html}}}

\Author{Evgeny FEIGIN~$^{\dag\ddag}$}

\AuthorNameForHeading{E.~Feigin}

\Address{$^\dag$~Mathematical Institute, University of Cologne, Weyertal 86-90, D-50931,
Cologne, Germany}

\Address{$^\ddag$~I.E.~Tamm Department of Theoretical Physics, Lebedev Physics Institute,\\
$\phantom{^\ddag}$~Leninski Prospect 53, Moscow,  119991, Russia}
\EmailD{\href{mailto:evgfeig@gmail.com}{evgfeig@gmail.com}}

\ArticleDates{Received July 04, 2008, in f\/inal form October
06, 2008; Published online October 14, 2008}

\Abstract{Let $L$ be the basic (level one vacuum) representation of the af\/f\/ine Kac--Moody Lie algebra
$\gh$. The $m$-th space $F_m$ of the PBW f\/iltration on $L$ is a linear span of vectors
of the form $x_1\cdots x_lv_0$, where $l\le m$, $x_i\in \gh$ and $v_0$ is a highest weight
vector of $L$. In this paper we give two descriptions of the associated graded space $L^{\rm gr}$ with respect
to the PBW f\/iltration. The ``top-down'' description deals with a structure of $L^{\rm gr}$ as
a representation of the abelianized algebra of generating operators. We prove that
the ideal of relations is generated by the coef\/f\/icients of the squared f\/ield
$e_\ta(z)^2$, which corresponds to the longest root $\ta$. The ``bottom-up'' description
deals with the structure of $L^{\rm gr}$ as a representation of the current algebra
$\g\T\C[t]$. We prove that each quotient $F_m/F_{m-1}$ can be f\/iltered by graded
deformations of the tensor products of $m$ copies of $\g$.}


\Classification{17B67}

\renewcommand{\thefootnote}{\arabic{footnote}}
\setcounter{footnote}{0}

\section{Introduction}
Let $\g$ be a f\/inite-dimensional simple Lie algebra, $\gh$ be the corresponding af\/f\/ine
Kac--Moody Lie algebra (see \cite{K1,Kum}). Let $L$ be the basic representation of $\gh$, i.e.\ an irreducible
level one module with a highest weight vector $v_0$ satisfying condition
$(\g\T\C[t])\cdot v_0=0$. The PBW f\/iltration $F_\bullet$ on the space $L$ is def\/ined as follows:
\[
F_0=\C v_0,\ F_{m+1}=F_m+\mathrm{span}\{x\cdot w:\ x\in\gh, \;w\in F_m\}.
\]
This f\/iltration was introduced in \cite{FFJMT} as a tool of study of vertex operators
acting on the space of Virasoro minimal models (see \cite{DMS}).
In this paper we study the associated graded space $L^{\rm gr}=F_0\oplus F_1/F_0\oplus\cdots$.
We describe the space $L^{\rm gr}$ from two dif\/ferent points of view: via ``top-down'' and
``bottom-up'' operators (the terminology of~\cite{KLMW}).

On one hand,
the space $L^{\rm gr}$ is a module over the Abelian Lie algebra $\g^{ab}\T t^{-1}\C[t^{-1}]$,
where $\g^{ab}$ is an Abelian Lie algebra whose underlying vector space is $\g$.
The module structure is induced from the action
of the algebra of generating ``top-down'' operators $\g\T t^{-1}\C[t^{-1}]$ on~$L$.
Thus $L^{\rm gr}$ can be identif\/ied with a polynomial ring on the space
$\g^{ab}\T t^{-1}\C[t^{-1}]$ modulo certain ideal. Our f\/irst goal is
to describe this ideal explicitly.

On the other hand,  all spaces $F_m$ are stable with respect to the
action of the subalgebra of annihilating operators $\g\T \C[t]$ (the ``bottom-up'' operators).
This gives
$\g\T\C[t]$ module structure on each quotient $F_m/F_{m-1}$. Our second goal
is to study these modules.
We brief\/ly formulate our results below.

For $x\in\g$ let $x(z)=\sum_{i<0} (x\T t^i) z^{-i}$ be a generating function
of the elements $x\T t^i\in\gh$, $i<0$. These series are also called f\/ields.
They play a crucial role in the theory of vertex operator algebras (see \cite{FK,K2, BF}). We will need the f\/ield which corresponds to the
highest root $\ta$ of $\g$. Namely, let $e_\ta\in\g$ be a highest weight vector
in the adjoint representation.
It is well known (see for instance \cite{BF}) that the coef\/f\/icients of $e_\ta(z)^2$ vanish on $L$.
It follows immediately that the same relation holds on $L^{\rm gr}$.
We note also that the Lie algebra $\g\simeq \g\T 1$
acts naturally on  $L^{\rm gr}$.
The following theorem is one of the central results of our paper.

\begin{theorem}\label{one}
$L^{\rm gr}$ is isomorphic to the quotient of the universal enveloping algebra
$\U(\g^{ab}\T t^{-1}\C[t^{-1}])$
by the ideal $I$, which is the minimal $\g\T 1$ invariant ideal containing all coefficients
of the series $e_\ta(z)^2$.
\end{theorem}

This proves the level one case of the conjecture from \cite{F2}. We note that for $\g=\slt$
this theorem was proved in \cite{FFJMT}. The generalization of this theorem for higher levels
and $\g=\slt$ is conjectured in \cite{F2}.

In order to prove this theorem and to make a connection with the
``bottom-up'' description we study the intersection of the PBW f\/iltration
with certain Demazure modules inside $L$. Recall (see \cite{D}) that by def\/inition a Demazure module
$D(\la)\hk L$ is
generated by extremal vector of the weight $\la$ with an action of the universal enveloping algebra
of the Borel subalgebra of $\gh$. We will only need the Demazure modules $D(N\ta)$.
These modules are invariant with respect to the current algebra $\g\T\C[t]$ and
provide a f\/iltration on $L$ by f\/inite-dimensional spaces:
$D(\ta)\hk D(2\ta)\hk\cdots =L$ (see \cite{FoL1, FoL2}; some special cases
are also contained in \cite{FF2,Ked}).
Let $F_m(N)=D(N\ta)\cap F_m$ be an intersection of the Demazure
module with the $m$-th space of the PBW f\/iltration. This gives a f\/iltration on $D(N\ta)$.
In order to describe the f\/iltration $F_\bullet(N)$ we
use a notion of the fusion product of $\g\T\C[t]$ modules (see \cite{FL,FF1}) and
the Fourier--Littelmann results \cite{FoL2}.

We recall that there exist two versions of the fusion procedure for modules over
the current algebras. The f\/irst version constructs a graded $\g\T\C[t]$ module
$V_1*\dots *V_N$
starting from the tensor product of cyclic $\g\T\C[t]$ modules $V_i$. The other version
also produces a graded $\g\T\C[t]$ module $V_1 **\cdots ** V_N$, but in this case all
$V_i$ are cyclic $\g$ modules. (We note that second version is a special case of the f\/irst one).
The fusion modules provide a useful tool for the study of the representation theory of current and
af\/f\/ine algebras (see \cite{AK, AKS,CLe, FKL,AKS,Ked,F1, F2,FJKLM, FoL2}).
In particular, Fourier and Littelmann proved
that there exists
an isomorphism of $\g\T\C[t]$ modules
$D(N\ta)\simeq D(\ta)*\cdots * D(\ta)$ ($N$ times).
Using this theorem and the $**$-version of the fusion procedure, we endow the space $D(N\ta)$
with a structure of the representation of the toroidal  current algebra
$\g\T\C[t,u]$ (see \cite{H, L} and references therein for some details on the
representation theory of the toroidal algebras).
This allows to prove our second main theorem:
\begin{theorem}\label{two}
The $\g\T\C[t]$ module $F_m(N)/F_{m-1}(N)$ can be filtered by
$\binom{N}{m}$ copies of the \linebreak \mbox{$m$-th}   fusion power of the adjoint representation of
$\g$. In particular,
$\dim F_m(N)/F_{m-1}(N)=\binom{N}{m} (\dim\g)^m.$
\end{theorem}

The paper is organized as follows.
In Section~\ref{sec1} we give the def\/inition of the PBW f\/iltration and of the induced
f\/iltration on Demazure modules.
In Section~\ref{sec2} we study tensor products of cyclic $\g\T\C[t]$ modules
endowed with a structure of representations of toroidal algebra. In particular,
we show that fusion product $D(1)^{**N}$ is well def\/ined.
In Section~\ref{sec3} the results of Section~\ref{sec2} are applied to the module $D(N)$. We prove
a graded version of the inequality  $\dim F_m(N)/F_{m-1}(N)\ge \binom{N}{m} (\dim\g)^m.$
In Section~\ref{sec4} the functional realization of the dual space $(L^{\rm gr})^*$ is given.
In Section~\ref{main} we combine all results of the previous sections and prove
Theorems~\ref{one} and~\ref{two}.
We f\/inish the paper with a list of the main notations.

\section[The PBW filtration]{The PBW f\/iltration}\label{sec1}

In this section we recall the def\/inition and basic properties  of the PBW f\/iltration
(see \cite{FFJMT}).

Let $\g$ be a simple Lie algebra, $\gh$ be the corresponding af\/f\/ine Kac--Moody
Lie algebra:
\[
\gh=\g\T\C[t,t^{-1}]\oplus \C K\oplus\C d.
\]
Here $K$ is a central element, $d$ is a degree element ($[d,x\T t^i]=-ix\T t^i$)
and
\[
[x\T t^i, y\T t^j]=[x,y]\T t^{i+j}+ i\delta_{i+j,0} (x,y)K,
\]
$x,y\in\g$, $(\cdot,\cdot)$ is a Killing form.
Let $L$ be the basic representation of the af\/f\/ine Lie algebra, i.e.\ level one
highest weight irreducible module with a highest weight vector $v_0$ satisfying
\[
(\g\T\C[t])\cdot v_0=0,\qquad Kv_0=v_0, \qquad dv_0=0,\qquad \U\big(\g\T t^{-1}\C[t^{-1}]\big)\cdot v_0=L.
\]
The operator $d\in\gh$  def\/ines a graded character of any subspace $V\hk L$ by the formula
\[
\ch_q V=\sum_{n\ge 0} q^n\dim\{v\in V:\ dv=nv\}.
\]
For $x\in\g$ we introduce a generating function
$x(z)=\sum_{i>0} (x\T t^{-i})z^i$ of the elements
$x\T t^i$, $i<0$.
We will mainly deal with the function $e_\ta(z)$, where $\ta$ is the
highest weight of $\g$ and $e_\ta\in\g$ is a highest weight element.
All coef\/f\/icients
\[
\sum_{\at{i+j=n}{i,j\le -1}} (e_\ta\T t^i)(e_\ta\T t^j)
\]
of the square of the series $e_\ta(z)$ are known to vanish on $L$
(this follows from the vertex operator realization of $L$ \cite[Theorem~$A$]{FK}).
Equivalently, $e_\ta(z)^2=0$ on $L$.

In what follows we will need a certain embedding of the basic $\slth$ module into $L$.
Namely, let~$\slt^\ta$ be a Lie algebra spanned by a $\slt$-triple $e_\ta$, $f_\ta$ and
$h_\ta$, where $e_\ta$ and $f_\ta$ are highest and lowest weight vectors in
the adjoint representation of $\g$. Then the restriction map def\/ines a~structure
of $\widehat{\slt^\ta}$ module on~$L$. In particular, the space $\U(\widehat{\slt^\ta})\cdot v_0$
is isomorphic to the basic representation of $\slth$, since the def\/ining relations
$(e_\ta\T t^{-1})^2v_0=0$ and $(\slt^\ta\T 1)v_0=0$ are satisf\/ied (see \cite[Lemma~2.1.7]{Kum}).

We now def\/ine the PBW f\/iltration $F_\bullet$ on $L$. Namely, let
\[
F_0=\C v_0,\qquad F_{m+1}=F_m+\mathrm{span}\{(x\T t^{-i}) w,\ x\in\g,\;
i>0,\; w\in F_m\}.
\]
Then $F_\bullet$ is an increasing f\/iltration
converging to $L$. We denote the associated graded space by~$L^{\rm gr}$:
\[
L^{\rm gr}=\bigoplus_{m\ge 0} L^{\rm gr}_m,\qquad
L^{\rm gr}_m=\gr_m F_\bullet=F_m/F_{m-1}.
\]

In what follows we denote by $\g^{ab}$ an Abelian Lie algebra with the underlying vector
space isomorphic to $\g$. We endow $\g^{ab}$ with a structure of $\g$ module via the
adjoint action of $\g$ on itself.
\begin{lemma}
The  action of $\g\T\C[t]$ on $L$ induces an action of the same algebra
on $L^{\rm gr}$. The  action of $\g\T t^{-1}\C[t^{-1}]$ on $L$ induces an action of the algebra
$\g^{ab}\T t^{-1}\C[t^{-1}]$ on $L^{\rm gr}$.
\end{lemma}
\begin{proof}
All  spaces $F_m$ are invariant with respect to the action of
$\g\T\C[t]$, since the condition  $(\g\T\C[t])\cdot v_0=0$ is satisf\/ied. Hence we obtain an
induced action on the quotient spaces $F_m/F_{m-1}$.

Operators $x\T t^i$, $i<0$ do not preserve $F_m$ but map it to $F_{m+1}$. Therefore, each
element $x\T t^i$, $i<0$ produce an operator acting from $L^{\rm gr}_m$ to $L^{\rm gr}_{m+1}$.
An important property of these operators on $L^{\rm gr}$ is that they mutually commute,
since the composition $(x\T t^i) (y\T t^j)$ acts from~$F_m$ to~$F_{m+2}$ but the commutator
$[x\T t^i, y\T t^j]=[x,y]\T t^{i+j}$ maps $F_m$ to $F_{m+1}$. Lemma is proved.
\end{proof}

The goal of our paper is to describe the structure of $L^{\rm gr}$ as a representation of $\g\T\C[t]$
and of $\g^{ab}\T t^{-1}\C[t^{-1}]$. It turns out that these two structures are closely related.

\begin{lemma}\label{fta}
Let $I\hk \U(\g^{ab}\T t^{-1}\C[t^{-1}])$ be the minimal $\g$-invariant ideal containing all
coefficients of the series $e_\ta(z)^2$. Then there exists a surjective homomorphism
\begin{gather}\label{*}
\U\big(\g^{ab}\T t^{-1}\C[t^{-1}]\big)/I\to L^{\rm gr}
\end{gather}
of $\g^{ab}\T t^{-1}\C[t^{-1}]$ modules mapping $1$ to $v_0$.
\end{lemma}
\begin{proof}
Follows from the relation $e_\ta(z)^2v_0=0$.
\end{proof}

One of our goals is to prove that the homomorphism (\ref{*}) is an isomorphism.

Recall (see \cite{FoL1}) that $L$ is f\/iltered by f\/inite-dimensional Demazure modules \cite{D}.
A Demazure module $D$ is generated by an extremal vector with the action of the
algebra of generating operators. We will only need special class of Demazure modules.
Namely, for $N\ge 0$ let $v_N\in L$ be the vector of weight $N\ta$ def\/ined by
\[
v_N=\big(e_\ta\T t^{-N}\big)^N v_0.
\]
We recall that $N\ta$ is an extremal weight for $L$ and thus $v_N$ spans
weight $N\ta$ subspace of $L$.
Let $D(N)\hk L$ be
the Demazure module generated by vector $v_N$.
Thus $D(N)$ is  cyclic $\g\T\C[t]$ module with  cyclic vector $v_N$.
It is known (see \cite{Kum,FoL1}) that these modules are embedded
into each other and the limit coincides with $L$:
\[
D(1)\hk D(2)\hk\cdots =L.
\]
We introduce an induced PBW f\/iltration on $D(N)$:
\begin{gather}\label{Ds} F_0(N)\hk F_1(N)\hk\cdots =D(N),\qquad
F_m(N)=D(N)\cap F_m.
\end{gather}
Obviously, each space $F_m(N)$ is $\g\T\C[t]$ invariant.

\begin{lemma}
$F_N(N)=D(N)$, but $F_{N-1}(N)\ne D(N)$.
\end{lemma}
\begin{proof}
First equality holds since $(e_\ta\T t^{-N})^N v_0=v_N$
and $D(N)$ is cyclic. In addition
$v_N\notin F_{N-1}$ because all weights of $F_m$ (as a representation of $\g\simeq \g\T 1$) are less
than or equal to $m\ta$.
\end{proof}

Consider the associated graded object
\[
\gr F(N)=\bigoplus_{m=0}^\infty \gr_m F(N),\qquad \gr_m F(N)=F_m(N)/F_{m-1}(N).
\]
We note that each space $\gr_m F(N)$ has a natural structure of $\g\T\C[t]$ module.

\section[$t^N$-filtration]{$\boldsymbol{t^N}$-f\/iltration} \label{sec2}

In this section we describe the f\/iltration (\ref{Ds}) using the generalization of the fusion product of $\g\T\C[t]$ modules from \cite{FL} and a theorem of
\cite{FoL2}. We f\/irst recall the def\/inition of the fusion product of $\g\T\C[t]$ modules.

Let $V$ be a $\g\T\C[t]$ module, $c$ be a complex number. We denote by
$V(c)$ a $\g\T\C[t]$ module which coincides with $V$ as a vector space and
the action is twisted by the Lie algebra homomorphism
\[
\phi(c): \ \g\T\C[t]\to\g\T\C[t],\quad x\T t^k\mapsto x\T (t+c)^k.
\]
Let $V_1,\dots, V_N$ be cyclic representations of the current algebra $\g\T\C[t]$ with
cyclic vectors $v_1,\dots,v_N$. Let $c_1,\dots,c_N$ be a set of pairwise distinct
complex numbers.
The fusion product $V_1(c_1)*\cdots * V_N(c_N)$ is a graded deformation of the usual tensor product
$V_1(c_1)\T \cdots\T V_N(c_N).$
More precisely,
let $\U(\g\T\C[t])_s$ be a natural grading on the universal enveloping
algebra coming from the counting of the $t$ degree. Because of the condition
$c_i\ne c_j$ for $i\ne j$, the tensor product $\bigotimes_{i=1}^N V_i(c_i)$
is a cyclic $\U(\g\T\C[t])$ module with a cyclic vector $\otimes_{i=1}^N v_i$. Therefore, the
grading on $\U(\g\T\C[t])$ induces an increasing fusion f\/iltration
\begin{gather}\label{fusfiltr}
\U(\g\T\C[t])_{\le s}\cdot (v_1\T\cdots\T v_N)
\end{gather}
on the tensor product.

\begin{definition}\label{D1}
The fusion product $V_1(c_1)*\cdots *V_N(c_N)$ of $\g\T\C[t]$ modules $V_i$  is an associated
graded $\g\T\C[t]$ module
with respect to the f\/iltration (\ref{fusfiltr}) on the tensor product
$V_1(c_1)\T\cdots\T V_N(c_N)$. We denote the $m$-th graded component by
$\gr_m (V_1(c_1)*\cdots * V_N(c_N)).$
\end{definition}

We note that in many special cases the $\g\T\C[t]$ module
structure of
the fusion product does not depend on the parameters $c_i$
(see for example \cite{AK, CL,FF1,FoL2,FKL}).
We then omit the parameters $c_i$ and denote the corresponding fusion
product simply by $V_1*\cdots *V_N$.

In what follows we will need a special but important case of the procedure described above.
Namely, let $V_i$ be cyclic $\g$ modules. One can
extend the $\g$ module structure to the $\g\T\C[t]$ module structure
by letting $\g\T t\C[t]$ to act by zero.
We denote the corresponding $\g\T\C[t]$ modules by $\overline{V}_i$.

\begin{definition}\label{D2}
Let $V_1,\dots,V_N$ be a set of cyclic $\g$ modules.
Then a $\g\T\C[t]$ module $V_1(c_1)**\cdots ** V_N(c_N)$
is def\/ined by the formula:
\[
V_1(c_1)**\cdots ** V_N(c_N)=\overline{V}_1(c_1)*\cdots * \overline{V}_N(c_N).
\]
We denote the $m$-th graded component by
$\gr_m (V_1(c_1)**\cdots ** V_N(c_N)).$
\end{definition}

\begin{remark}
The fusion procedure described in Def\/inition \ref{D2} can be reformulated as follows.
One starts with a tensor product of evaluation $\g\T\C[t]$ modules
$V_i(c_i)$, where $x\T t^k$ acts on $V_i$ by $c_i^kx$ (we evaluate
$t$ at the point $c_i$). Then one constructs the fusion f\/iltration and
associated graded module (according to Def\/inition~\ref{D1}).
\end{remark}

\begin{remark}
Let us stress the main dif\/ference between Def\/initions~\ref{D1} and~\ref{D2}.
Def\/inition~\ref{D1} gets as an input a set of cyclic representations of
the current algebra $\g\T\C[t]$ and as a result produces a $t$-graded $\g\T\C[t]$
module.
The input of Def\/inition~\ref{D2} is a set of cyclic $\g$ modules and an output
is again a $t$-graded $\g\T\C[t]$ module.
\end{remark}

We now recall a theorem from \cite{FoL2} which uses the fusion procedure to construct
the Demazure module $D(N)$ starting from $D(1)$.

\begin{theorem}[\cite{FoL2}] \label{FL}
The $N$-th fusion power $D(1)^{*N}$ is independent on the evaluation
parame\-ters~$c_i$. The $\g\T\C[t]$ modules $D(1)^{*N}$ and $D(N)$ are isomorphic.
\end{theorem}

We recall that as a $\g$ module  $D(1)$ is isomorphic to the direct sum of trivial and
adjoint representations, $D(1)\simeq \C\oplus\g$. The trivial representation
is annihilated by $\g\T\C[t]$, the adjoint representation is annihilated by
$\g\T t^2\C[t]$ and $\g\T\C[t]$ maps $\g$ to $\C$.

Our idea is to combine the theorem above and Def\/inition~\ref{D2} with
$\g$ being the current algebra $\g\T\C[u]$ (see also \cite{CLe}, where Def\/inition \ref{D2}
is used in af\/f\/ine settings). Def\/inition~\ref{D2}
works for arbitrary $\g$ and produces a representation of an algebra with an
additional current variable. In particular, starting from the $\g\T \C[u]$ modules
$V_i=D(1)$ and an $N$ tuple of pairwise distinct complex numbers $c_1,\dots,c_N$
one gets a new bi-graded $\g\T\C[t,u]$ module.
The resulting module can be obtained
from the Demazure module $D(N)$ by a rather simple
procedure which we are going to describe now.

Let $V_1,\dots, V_N$ be cyclic representations of the algebra $\g\T\C[t]/\bra t^2\ket$.
Hence $V_1(c_1)*\cdots *V_N(c_N)$ is a cyclic $\g\T\C[t]/\langle t^{2N}\rangle$ module.
We consider a decreasing f\/iltration $\U(\g\T\C[t])^j$ on the universal enveloping
algebra def\/ined by
\begin{gather}\label{Gdef}
\U(\g\T\C[t])^0=\U(\g\T\C[t]),\qquad \U(\g\T\C[t])^{j+1}=\big(\g\T t^N\C[t]\big)\U(\g\T\C[t])^j.
\end{gather}
This f\/iltration induces a decreasing f\/iltration $G^j$
on the fusion product $V_1(c_1)*\dots *V_N(c_N)$ (since it is cyclic
$\U(\g\T\C[t])$ module). $G^\bullet$ will be also referred to as a $t^N$-f\/iltration.
Consider the associated graded space
\begin{gather}\label{Gad}
\gr G^\bullet =\bigoplus_{j=0}^N \gr^j G^\bullet,\qquad \gr^j G^\bullet=G^j/G^{j+1}.
\end{gather}
Since each space $G^j$ is $\g\T\C[t]$ invariant one gets a structure
of $\g\T \C[t]/\bra t^N\ket$ module on each space $G^j/G^{j+1}$.  In addition
an element from $\g\T t^N\C[t]/\langle t^{2N}\rangle$ produces a degree $1$  operator on~(\ref{Gad})
mapping $\gr^j G^\bullet$ to $\gr^{j+1} G^\bullet$.
We thus obtain a structure of $\g\T\C[t,u]/\bra t^N,u^2\ket$ module on~(\ref{Gad}),
where $\g\T u\C[t]$ denotes an algebra of degree one operators on $\gr G^\bullet$
coming from the action of $\g\T t^N\C[t]/\langle t^{2N}\rangle$.

On the other hand let us consider the modules $V_i$ as representations of the Lie
algebra $\g\T\C[u]/\bra u^2\ket$ (simply replacing $t$ by $u$). We denote these modules
as $V_i^u$. Then the bi-graded
tensor product $V^u_1(c_1)**\cdots ** V^u_N(c_N)$ is a representation of the Lie algebra
$\g\T\C[t,u]/\bra t^N, u^2\ket$.

\begin{proposition}\label{t^N}
We have an isomorphism of $\g\T\C[t,u]/\langle t^N,u^2\rangle$ modules
\begin{gather}\label{**}
\gr G^\bullet\simeq V^u_1(c_1)**\cdots **V^u_N(c_N).
\end{gather}
\end{proposition}
\begin{proof}
The idea of the proof is as follows.
We start with the tensor product $V_1\T\cdots\T V_N$ and apply the fusion f\/iltration.
Afterwards we apply the $t^N$-f\/iltration $G^\bullet$. Combining these operations
with certain changes of basis of current algebra we arrive at the def\/inition of
the bi-graded module $V^u_1**\cdots ** V^u_N$. We give details below.

For an element $x\T t^i\in \g\T\C[t]$ let $(x\T t^i)^{(j)}$ be the operator on
the tensor product $V_1\T\cdots\T V_N$ def\/ined by
\[
(x\T t^i)^{(j)}= \Id^{\T j-1} \T (x\T t^i)\T\Id^{\T N-j} ,
\]
i.e.\ $(x\T t^i)^{(j)}$ acts as $x\T t^i$ on $V_j$ and as an identity operator
on the other factors. In order to construct  the fusion product one starts with
the  operators
\begin{gather}
\label{1fus}
A(x\T t^i)=\sum_{j=1}^N \big(x\T (t+c_j)^i\big)^{(j)},
\end{gather}
where $c_j$ are pairwise distinct numbers. These operators def\/ine an action of the
algebra $\g\T\C[t]$ on the tensor product $V_1(c_1)\T\cdots \T V_N(c_N)$.
Since $x\T t^i$ with $i>1$ vanish on $V_j$ we obtain
\[
A(x\T t^i)=\sum_{j=1}^N c_j^i (x\T 1)^{(j)}+ \sum_{j=1}^N ic_j^{i-1} (x\T t)^{(j)}.
\]
The next step is to pass to  the associated graded module with respect to
the fusion f\/iltration.
By def\/inition, operators of the form
\begin{gather}\label{lc}
A\big(x\T t^k\big)+\text{ linear combination of } A\big(x\T t^l\big), \qquad l<k
\end{gather}
act on the associated graded module in the same way as $A(x\T t^k)$
(the lower degree term vanish after passing to the associated graded space).
We are going to f\/ix special linear changes in (\ref{lc}) for $N\le k <2N$
which make the expressions for $A(x\T t^k)$ simpler.

Let $\al_0,\al_1,\dots,\al_{N-1}$ be numbers such that for all $1\le j\le N$
\[
c_j^N+\sum_{i=0}^{N-1} \al_i c_j^i=0.
\]
We state that
\begin{gather}\label{1}
A\big(x\T t^{N+s}\big)+\sum_{i=0}^{N-1} \al_i A\big(x\T t^{s+i}\big)= \sum_{j=1}^N c_j^s (x\T t)^{(j)}\prod_{k\ne j} (c_j-c_k)
\end{gather}
for all $0\le s\le N-1$.
Let
\[
f(x)=x^N+\al_{N-1}x^{N-1}+\dots + \al_0.
\]
Then $f(x)=\prod_{k=1}^N (x-c_k)$. Therefore, for the derivative $c_j^sf'(c_j)$ one gets
\[
c_j^sf'(c_j)=Nc_j^{N+s-1}+\sum_{i=0}^{N-1} i\al_ic_j^{i+s-1}=c_j^s\prod_{k\ne j} (c_j-c_k).
\]
This proves (\ref{1}).

Using  formula (\ref{1}), we replace operators $A(x\T t^i)$,
$0\le i< 2N$ by operators $B(x\T t^i)$ as follows
\begin{gather}
\label{u0}
B\big(x\T t^i\big)=\sum_{j=1}^N c_j^i (x\T 1)^{(j)} + \sum_{j=1}^N ic_j^{i-1} (x\T t)^{(j)},\qquad 0\le i< N,\\
\nonumber
B\big(x\T t^{N+i}\big)=\sum_{i=1}^N c_i^s (x\T t)^{(i)} \prod_{k\ne j} (c_j-c_k),\qquad 0\le i<N,
\end{gather}
thus performing the linear change (\ref{lc}).
So we redef\/ine half of the operators $A(x\T t^i)$ and leave the other half unchanged.

In order to construct the left hand side of (\ref{**}) one f\/irst applies the fusion
f\/iltration to the algebra of operators $B(x\T t^i)$ and afterwards proceeds with
the $t^N$-f\/iltration.
The last step means that the subtraction of a linear combination of
the operators $B(x\T t^{N+i})$, $0\le i< N$ from $B(x\T t^i)$, $0\le i<N$ does not
change the structure of the resulting module. Redef\/ining the operators (\ref{u0})
we arrive at the following operators:
\begin{gather*}
C\big(x\T t^i\big)=\sum_{j=1}^N c_j^i (x\T 1)^{(j)},\qquad  0\le i< N,\\
C\big(x\T t^{N+i}\big)=\sum_{i=j}^N c_j^i (x\T t)^{(j)},\qquad 0\le i<N.
\end{gather*}
Note that we can remove constants $\prod_{k\ne j} (c_j-c_k)$ from
$B(x\T t^{N+i})$ since this procedure corresponds simply to
rescaling the variable $t$ in each $V_i$.

Summarizing all the formulas above we arrive at the following two steps construction of the
left hand side of (\ref{**}):
\begin{itemize}\itemsep=0pt
\item apply the fusion procedure to the operators $C(x\T t^i)$, $0\le i<2N$,
\item attach a $u$-degree $1$ to each of the operators $C(x\T t^{N+i})$, $0\le i< 0$.
\end{itemize}
In order to construct the right hand side of (\ref{**}) one uses the same procedure
with $C(x\T t^{N+i})$ being operators which correspond to $x\T ut^i$ (see (\ref{Gdef}) and
(\ref{Gad})).
Thus we have shown that the associated graded to the fusion product with respect to the
$t^N$-f\/iltration is isomorphic to the module $V_1**\dots ** V_N$.
Proposition is proved.
\end{proof}

\begin{corollary}
The fusion product $D(1)**\cdots ** D(1)$ does not depend on the evaluation para\-me\-ters.
\end{corollary}

\begin{corollary}
The fusion product of the adjoint representations
\[
\g(c_1)**\cdots **\g(c_N)
\]
is independent of the parameters $c_1,\dots,c_N$.
\end{corollary}
\begin{proof}
By def\/inition, the zeroth graded component with respect to the $t$ grading
of the module
\[
D(1)(c_1)**\cdots ** D(1)(c_N)
\]
is isomorphic
to the fusion product $\g(c_1)**\dots **\g(c_N)$. From the Proposition above we obtain an
isomorphism
\[
\g(c_1)**\dots **\g(c_n)\simeq D(N)/ (\g\T t^N\C[t]) D(N)
\]
for any $c_1,\dots,c_n$. Thus the left hand side is independent of $c_i$.
\end{proof}

We f\/inish this section  introducing an ``energy'' operator~$\bar d$ on the fusion product.
The opera\-tor~$\bar d$ acts by a constant $m$
on the graded component
\[
\gr_m (V_1(c_1) *\dots * V_N(c_N)).
\]
The operator $\bar d$ def\/ines a graded character of $V_1(c_1)*\dots *V_n(c_n)$ by the standard formula formula
\[
\ov{\ch}_q V_1(c_1)*\dots *V_n(c_n)=\sum_{n\ge 0} q^n\dim\{v:\ \bar dv=nv\}.
\]
An analogous formula def\/ines a character of $V_1(c_1)**\dots **V_n(c_n)$.
We note that this grading has nothing to do with the Cartan grading coming from the
action of the Cartan subalgebra. We do not consider the latter grading in this paper.

\begin{remark}\label{ch}
Let $V_i\simeq D(1)$ for all $i$. Then the fusion module
is independent on the evaluation parameters and $D(1)^{*N}$ is isomorphic to the
Demazure module $D(N)$ (see \cite{FoL2}). By the very def\/inition we have an embedding $D(N)\hk L$. Thus
both operators $d$ and $\bar d$ are acting on $D(N)$, satisfying the relations
\[
[d,x\T t^i]=-ix\T t^i,\qquad [\bar d,x\T t^i]=ix\T t^i.
\]
Since $dv_N=N^2v_N$ and $\bar dv_N=0$ we have a
simple identity $\bar d=N^2-d$.
\end{remark}

\section{Demazure modules}\label{sec3}
In this section we study the fusion f\/iltration on the tensor product $D(1)^{\T N}$
and the induced PBW f\/iltration on the Demazure modules $D(N)$. We also derive some
connections between these f\/iltrations.

Let $D^u(1)$ be the $\g\T\C[u]/\bra u^2\ket$ module obtained from $D(1)$ by
substituting $u$ instead of $t$. In particular, $(D^u(1))^{**N}$ is a $(t,u)$ bi-graded
representation of the algebra $\g\T\C[t,u]/\bra t^N,u^2\ket$.
Here $t$-grading is exactly the fusion grading $\gr_j D^u(1)*\dots * D^u(1)$
and the $u$-grading comes from the grading on $D^u(1)$, which assigns degree zero to $\g$
and degree one to $\C v_0$.
We consider  the decomposition
\[
(D^u(1))^{**N}=\bigoplus_{m=0}^N W(m)
\]
into the graded components with respect to the $u$-grading.
(The $u$-grading is bounded from above by $N$ since the $u$-grading
in each of $D(1)$ could be either $0$ or $1$).
Note that each space~$W(m)$ is a representation of $\g\T\C[t]$.
We want to show that~$W(m)$ can be f\/iltered by certain number
of copies of the fusion product $\g^{** N-m}$. The precise statement is given in the
following proposition.

\begin{proposition}\label{g*s}
Let $D^u(1)^{**N}$ be a bi-graded tensor product of $N$ copies of  $\g\T\C[u]$-module~$D(1)$.  Then
\begin{itemize}\itemsep=0pt
\item
For any $0\le m\le N$ the $\g\T\C[t]$ module $W(m)$
can be filtered by $\binom{N}{m}$ copies of the $\g\T\C[t]$-module $\g^{**N-m}$.
\item  The cyclic vectors of the modules $\g^{**N-m}$
above are the images of the vectors
\begin{gather}\label{u}
\big(f^\ta\T ut^{i_1}\big)\cdots \big(f^\ta\T ut^{i_m}\big) v_N,\qquad 0\le i_1\le\dots\le i_m\le N-m.
\end{gather}
\end{itemize}
\end{proposition}

\begin{remark}
We f\/irst give a non rigorous, but conceptual explanation of the statement of the
proposition above. Recall that $D^u(1)$ is isomorphic to $\g\oplus\C$
as a $\g$ module. Let $v_1,v_0\in D(1)$ be highest weight vectors of $\g$ and $\C$
respectively. Then $(f_\ta\T u)v_1=v_0$ and $(f_\ta\T u)^2v_0=0$. Therefore,
after making the fusion $D^u(1)^{**N}$, the tensor product of $N$ copies of
$2$-dimensional vector space $\mathrm{span}\{v_0,v_1\}$ will be deformed into
the $N$-fold fusion product of two-dimensional representations of the algebra
$\C[f_\ta]$. The set of vectors (\ref{u}) represents a basis of this
fusion product (see \cite{CP,FF1}). Hence $D^u(1)^{*N}$ is equal
to the $\U(\g\T\C[t])$ span of the vectors of the form~(\ref{u}).
We now want to describe the space
$\U(\g\T\C[t])\cdot w_m$, where $w_m$ is of the form~(\ref{u}) with
exactly $m$ factors. Note that $w_m$ is a linear combination of the vectors of the form
\[
v_{i_1}\T\cdots \T v_{i_N},
\]
where $i_\al$ equal $0$ or $1$ and the number of $\al$ such that $i_\al=0$ is equal
to $m$. This means that the space
 $\U(\g\T\C[t])\cdot w_m$ is embedded into the direct sum of $\binom{N}{m}$ copies
 of the tensor product $\g^{\T m}$. Hence it is natural to expect that after
 passing to the associated graded object with respect to the fusion f\/iltration
one arrives at the fusion product $\g^{**m}$. This is not exactly true.
In order to make the statement precise one additional f\/iltration is needed (that is the reason
why $W(m)$ is not the direct sum of the fusion products, but rather can be f\/iltered by
these modules).
\end{remark}

We now give the proof of Proposition \ref{g*s}.
\begin{proof}
As a starting point we note the isomorphism of $\g\T\C[t]$-modules
\[
W(0)\simeq \g^{**N}.
\]
In fact,
$D(1)\simeq \g\oplus\C v_0$ with cyclic vector being the highest weight vector of $\g$.
The algebra $\g\T u$ maps $\C$ to zero and $\g$ to $\C$.
In particular,
\[
(f_\ta\T u)\cdot v_1=v_0.
\]
Hence if we do not apply
operators with positive powers of $u$ (i.e. we consider the space $W(0)$)
we arrive at the usual fusion product of $N$ copies of $\g$.

We now introduce a decreasing f\/iltration $W^j(m)$ such that the associated graded object
is isomorphic to the direct sum of $\binom{N}{m}$ copies of $\g^{**m}$.
Let
\[
w_{i_1,\dots,i_m}= \big(f^\ta\T ut^{i_1}\big)\dots \big(f^\ta\T ut^{i_m}\big) v_N,\qquad 0\le i_1\le\dots\le i_m\le N-m.
\]
We set
\[
W^n(m)=\U(\g\T\C[t])\cdot\mathrm{span}\{ w_{i_1,\dots,i_m}:\ i_1+\cdots +i_m\ge n\}.
\]
In particular, $W^0(m)=W(m)$ and each space $W^n(m)$ is $\g\T\C[t]$ invariant.
We state that the associated graded space
\begin{gather}\label{W(s)}
W^0(m)/W^1(m)\oplus W^1(m)/W^2(m)\oplus\cdots
\end{gather}
is isomorphic to the direct sum of $\binom{N}{m}$ copies of the modules $\g^{**N-m}$.
Moreover the highest weight vectors of these modules are exactly the images of
$w_{i_1,\dots,i_m}$.
We prove this statement for $m=1$. The proof for other $m$ is very similar.

For $m=1$ we have
\[
w_i=f_\ta\T ut^i=\sum_{j=1}^N c_j^i v_1^{\T j-1}\T v_0\T v_1^{\T N-j},\qquad i=0,\dots,N-1.
\]
We want to show that
\[
W^i(1)/W^{i+1}(1)\simeq \g^{**N-1},\qquad W^i(1)=\U(\g\T\C[t])\cdot \mathrm{span}\{w_i,\dots,w_{N-1}\}.
\]
Let $\al_{i,j}$, $1\le i,j\le N$ be some numbers. Denote
\[
\tilde W^i(1)=\U(\g\T\C[t])\cdot \mathrm{span}\left\{\sum_{j=1}^N \al_{i+1,j} w_i,\dots, \sum_{j=1}^N \al_{N,j} w_{N-1}\right\}.
\]
We state that
\begin{gather}\label{tilde}
W^i(1)/W^{i+1}(1)\simeq \tilde W^i(1)/\tilde W^{i+1}(1)
\end{gather}
as $\g\T\C[t]$ modules. In fact, adding to $w_i$ a linear combination of $w_j$ with
$j<i$ does not change the $\g\T\C[t]$ module structure because of the fusion
f\/iltration. Because of the f\/iltration (\ref{W(s)}), this is still true if
one adds a linear combination of $w_j$ with $j>i$. Thus we conclude that we can replace
each vector $w_i$ by an arbitrary linear combination.
In particular, there exist numbers~$\al_{i,j}$ such that
\[
\sum_{j=1}^N \al_{i-1,j} w_i=v_0^{\T i} \T v_1\T v_0^{\T N-1-i}.
\]
Since $v_1$ is a highest weight vector of the adjoint representation of $\g$ and
$v_0$ spans the trivial representation, we arrive at the fact that
\[
\tilde W^i(1)/\tilde W^{i+1}(1)\simeq \g^{** N-1}.
\]
Because of the isomorphism (\ref{tilde}), the $m=1$ case of the proposition is proved.
\end{proof}

We are now going to connect the f\/iltrations $G^\bullet(N)$
and the induced PBW f\/iltration $F_\bullet(N)$ on the Demazure modules $D(N)$. We use
Proposition \ref{t^N} for $V_i=D(1)$.

\begin{lemma}\label{GF}
$G^m$ is a subspace of $F_{N-m}(N)$.
\end{lemma}
\begin{proof}
We f\/irst note that for the cyclic vector $v_N\in D(N)$ we have
\begin{gather}\label{v_N}
v_N=\big(e_\ta\T t^{-N}\big)^N v_0.
\end{gather}
Therefore, $G^0=F_N(N)=D(N)$ and our Lemma is true for $m=0$.

In general, we need to prove that
\begin{gather}\label{x1}
\big(x_1\T t^{N+i_1}\big)\cdots \big(x_m\T t^{N+i_m}\big)v\in F_{N-m}
\end{gather}
for any $v\in D(N)$ and $x_1,\dots, x_m\in\g$, $i_1,\dots,i_m\ge 0$.
Since $v\in D(N)$ there exists an element $r\in\U(\g\T\C[t])$ such that
$v=rv_N$. Because of (\ref{v_N}), the inclusion (\ref{x1}) follows
from the following statement:
\begin{gather}\label{FN-s}
\big(x_1\T t^{N+i_1}\big) \cdots \big(x_m\T t^{N+i_m}\big)\big(y_1\T t^{j_1}\big)\cdots \big(y_m\T t^{j_m}\big) \big(e_\ta\T t^{-N}\big)^N v_0\in F_{N-m}
\end{gather}
for arbitrary $x_\al,y_\be\in\g$ and $i_\al,j_\be, n\ge 0$.
Since $(\g\T\C[t])\cdot v_0=0$ the expression
\[
\big(y_1\T t^{j_1}\big)\cdots \big(y_n\T t^{j_n}\big) \big(e_\ta\T t^{-N}\big)^N v_0
\]
is equal to a linear combination of the vectors of the form
\[
\big(z_1\T t^{-N+l_1}\big)\cdots \big(z_N\T t^{-N+l_N}\big)v_0
\]
where $z_i\in\g$ and $l_i\ge 0$. One can easily see that because of the condition
$(\g\T\C[t])v_0=0$ the expression of the form
\[
\big(x\T t^{N+i}\big)\big(z_1\T t^{-N+l_1}\big)\cdots \big(z_N\T t^{-N+l_N}\big)v_0
\]
can be rewritten as a linear combination of the monomials of the form
\[
\big(z_1\T t^{-N+l_1}\big)\cdots \big(z_{N-1}\T t^{-N+1+l_{N-1}}\big)v_0.
\]
Iterating this procedure we arrive at (\ref{FN-s}). This proves (\ref{x1}) and
hence our lemma is proved.
\end{proof}

Recall a notation for $q$-binomial coef\/f\/icients
\[
\qbin{n}{m}=\frac{(q)_n}{(q)_m(q)_{n-m}},\qquad (q)_a=(1-q)\cdots (1-q^a).
\]
For two $q$ series $f(q)=\sum f_nq^n$ and $g(q)=\sum g_nq^n$ we write $f\ge g$ if
$f_n\ge g_n$ for all $n$.
\begin{corollary}
The following character inequality holds:
\begin{gather}\label{DsN}
\ch_q \gr_m F_\bullet(N)\ge q^{m^2}\qbin{N}{m}  \ov{\ch}_{q^{-1}} \g^{**m}.
\end{gather}
\end{corollary}
\begin{proof}
Because of the Lemma above, it suf\/f\/ices to show that
\[
\ch_q G^m/G^{m+1}=q^{(N-m)^2}\qbin{N}{m}  \ov{\ch}_{q^{-1}} \g^{**N-m}.
\]
Note that $dv_N=N^2v_N$. Therefore, the graded character of the space of cyclic
vectors $w_{i_1,\dots,i_m}$ from Proposition~\ref{g*s} (where the $q$-degree of
$u$ is f\/ixed to be equal to $N$ according to Proposition~\ref{t^N})  is given by
$q^{m^2}\qbin{N}{m}$. Multiplying by the $\bar d$ character of $\g^{**N-m}$ with respect to $\bar d$
(see Remark~\ref{ch}), we arrive at our Corollary.
\end{proof}

\begin{corollary}\label{Fs}
The following character inequality holds:
\begin{gather}\label{lim}
\ch_q \gr_m F_\bullet\ge q^{m^2}\frac{1}{(q)_m}  \ov{\ch}_{q^{-1}} \g^{**m}.
\end{gather}
\end{corollary}
\begin{proof}
Follows from limit formulas
\[
\lim_{N\to \infty} D(N)=L,\qquad \lim_{N\to\infty} \qbin{N}{m}=\frac{1}{(q)_m}.\tag*{\qed}
\]\renewcommand{\qed}{}
\end{proof}

\begin{remark}
In the next section we prove that~(\ref{lim}) is an equality.
In Section~\ref{main} we prove that~(\ref{DsN}) is also an equality.
\end{remark}

\section{Dual functional realization}\label{sec4}
We now consider the restricted dual space to the ``expected'' PBW f\/iltered space $L^{\rm gr}$.
Let
\[
\U^{ab}=\C[\g^{ab}\T t^{-1}\C[t^{-1}]]
\]
 be a space of polynomial functions on
the space $\g^{ab}\T t^{-1}\C[t^{-1}]$ (recall that $\g^{ab}$ is an Abelian Lie
algebra with the underlying vector space isomorphic to $\g$).
The algebra $\U^{ab}$ is an abelinization of the universal enveloping algebra
$\U(\g\T t^{-1}\C[t^{-1}])$ of the Lie algebra of generating operators (due to the
PBW theorem).
We note that $\g$ acts on the space $\g^{ab}\T t^{-1}\C[t^{-1}]$ via the adjoint
representation on $\g^{ab}$.

Let $I\in\U^{ab}$ be the minimal $\g$ invariant ideal, which contains all coef\/f\/icients
of $e^\ta(z)^2$, i.e. all elements of the form
\[
\sum_{1\le i\le n} \big(e_\ta\T t^{-i}\big)\big(e_\ta\T t^{-n-1+i}\big), \qquad n=1,2,\dots.
\]
Denote
\[
L'= \U^{ab}/I.
\]
This space can be decomposed according to the number of variables in a monomial:
\[
L'=\bigoplus_{m=0}^\infty L'_m=\bigoplus_{m=0}^\infty \mathrm{span}\big\{x_1\T t^{i_1}\cdots x_m\T t^{i_m},\;
x_i\in\g^{ab}\big\}.
\]
The operator $d\in\gh$ induces a degree operator on $L'$. We denote this operator by the same symbol.
There exists a surjective homomorphism of
$\g^{ab}\T t^{-1}\C[t^{-1}]$ modules $L'\to L^{\rm gr}$
(see Lemma~\ref{fta}).
Our goal is to show that $L'\simeq L^{\rm gr}$.
Let $(L'_m)^*$ be a restricted dual space:
\[
(L'_m)^*=\bigoplus_{n\ge 0} (L'_{m,n})^*,\qquad L'_{m,n}=\{v\in L'_m:\ dv=nv\}.
\]
We construct the functional realization of $L'_m$ using currents
\[
x(z)=\sum_{i>0} \big(x\T t^{-i}\big) z^i
\]
for $x\in\g$.
Following \cite{FFJMT} we consider a map
\[
\al_m: \ (L'_m)^*\to \C[z_1,\dots,z_m]\T \g^{\T m},\qquad \phi\mapsto r_\phi
\]
from the dual space $(L'_m)^*$ to the space of polynomials in $m$ variables with values
in the $m$-th tensor power of the space $\g$. This map is given by the formula
\[
\bra r_\phi, x_1\T\cdots\T x_m\ket=\phi(x_1(z_1)\cdots x_m(z_m)),
\]
where brackets in the left hand side denote
the product of non degenerate Killing forms on $n$ factors $\g$.
Our goal is to describe the image of $\al_m$. We f\/irst formulate
the conditions on $r_\phi$, which follow from the def\/inition of $L'$
and then prove that these conditions determine the image of $\al_m$.
We prepare some notations f\/irst.

Consider the decomposition of the tensor square $\g\T \g$ into the direct sum of $\g$ modules:
\begin{gather}\label{decomp}
\g\T\g=V_{2\ta}\oplus {\bigwedge}^2\g\oplus S^2\g/V_{2\ta}.
\end{gather}
Here $V_{2\ta}$ is a highest weight $\g$-module with a highest
weight $2\ta$ embedded into $S^2\g$ via the map
\[
V_{2\ta}\simeq \U(\g)\cdot (e_\ta\T e_\ta)\hk S^2\g.
\]

\begin{lemma}
For the  module $\g**\g$ we have
\[
\gr_0(\g**\g)=V_{2\ta},\qquad \gr_1 (\g**\g)={\bigwedge}^2\g,\qquad \gr_2 (\g**\g)=S^2\g/V_{2\ta}
\]
and all other graded components vanish.
\end{lemma}
\begin{proof}
We f\/irst show that $\gr_n (\g**\g)=0$ for $n>2$. Let $c_1$, $c_2$ be evaluation constants which
appear in Def\/inition~\ref{D2} of the fusion product.
Recall that $\g**\g$ is independent of the evaluation parameters. So we can set $c_1=1$, $c_2=0$.
Then the second space of the fusion f\/iltration is given~by
\begin{gather}\label{2}
\U(\g)\cdot\left(\mathrm{span}\{[x_1,[x_2,e_\ta]],\ x_1,x_2\in\g \}\T e_\ta\right),
\end{gather}
where $\U(\g)$ acts on the tensor product $\g\T\g$ diagonally.
But $[f_\ta,[f_\ta,e_\ta]]=-2f_\ta$ and hence
\[
\mathrm{span}\{[x_1,[x_2,e_\ta]],\ x_1,x_2\in\g \}=\g
\]
(since the left hand side is invariant with respect to the subalgebra of $\g$
of annihilating
operators and contains the lowest weight vector $f_\ta$ of the
adjoint representation). Therefore, (\ref{2}) coincides with $\g\T\g$.

We now compute three nontrivial graded components $\gr_0 (\g**\g)$, $\gr_1 (\g**\g)$ and $\gr_2 (\g**\g)$.
From the def\/inition of the fusion f\/iltration we have
\[
\gr_0 (\g**\g)= \U(\g)\cdot (e_\ta\T e_\ta)\simeq V(2\ta).
\]
We now redef\/ine the
evaluation parameters by setting $c_1=1$, $c_2=-1$. Then the formula
for the operators $x\T t$ acting on $\g\T\g$ is given by $x\T\Id- \Id\T x$. These operators
map $S^2(\g)$ to~${\bigwedge}^2\g$ and vice versa.We conclude that
\[
\gr_1 (\g**\g)\hk{\bigwedge}^2\g,\qquad \gr_2 (\g**\g) \hk S^2(\g).
\]
Now our Lemma follows from the equality $\gr_n (\g**\g)=0$ for $n>2$.
\end{proof}

\begin{lemma}
For any $\phi\in (L'_m)^*$ the image $r_\phi$ is divisible by the product $z_1\cdots z_m$ and
satisfies the vanishing
condition
\begin{gather}\label{van}
\bra r_\phi, V_{2\ta}\T \g^{\T m-2}\ket_{z_1=z_2}=0
\end{gather}
and the symmetry condition
\begin{gather}\label{sym}
\sigma r=r, \qquad \sigma\in S_m,
\end{gather}
where the permutation group $S_m$ acts simultaneously on the set of variables
$z_1,\dots,z_m$ and on the tensor product $\g^{\T m}$.
\end{lemma}
\begin{proof}
The product $z_1\cdots z_m$ comes from the condition that the highest weight vector
is annihilated by $\g\T \C[t]$, so for any $x\in\g$ the series
$x(z)$ starts with $z$.
The condition (\ref{van}) follows from the relation $e_\ta(z)^2=0$ and $\g$-invariance of the
ideal $I$. The symmetry condition follows from the commutativity of the algebra $\U^{ab}$.
\end{proof}

We denote by
\[
V_m\hk z_1\cdots z_m\C[z_1,\dots,z_m]\T\g^{\T m}
\]
the space of functions satisfying
conditions (\ref{van}) and (\ref{sym}).
In the following lemma we endow $V_m$ with structures of representation
of the ring of symmetric polynomials
\[
P_m^{\rm sym}=\C[z_1,\dots,z_m]^{S_m}
\]
and of the current algebra $\g\T\C[t]$.

\begin{lemma}\label{123}
There exists natural structures of representations of $P_m^{\rm sym}$ and of
$\g\T\C[t]$ on $V_m$ defined by the following rule:
\begin{itemize}\itemsep=0pt
\item $P_m^{\rm sym}$ acts on $V_m$ by multiplication on the first tensor factor.
\item Lie algebra $\g\T\C[t]$ acts on $V_m$ by the formula
\[
x\T t^k \quad \text{acts  as} \ \ \sum_{i=1}^n z_i^k\T x^{(i)}.
\]
\end{itemize}
The actions of $P_m^{\rm sym}$ and of $\g\T\C[t]$ commute.
\end{lemma}
\begin{proof}
A direct computation.
\end{proof}

Lemma~\ref{123} gives a structure of $\g\T\C[t]$ module
on the quotient space
$V_m/P^{\rm sym}_{m+}V_m$, where the subscript $+$ denotes the space of polynomials
of positive degree. We will show that the dual to this module is isomorphic to
$\g^{**m}$. We f\/irst consider the $m=2$ case.

\begin{lemma}\label{n=2}
We have an isomorphism of representations of $\g\T\C[t]$
\[
\left(V_2/P^{\rm sym}_+V_2\right)^*\simeq \g**\g.
\]
\end{lemma}
\begin{proof}
Let $r$ be an element of $V_2$. Using the decomposition (\ref{decomp}),
we write $r=r_0+r_1+r_2$, where
\begin{gather*}
r_0\in z_1z_2\C[z_1,z_2]\T V_{2\ta},\\
r_1\in z_1z_2\C[z_1,z_2]\T \Lambda^2\g,\\
r_2\in z_1z_2\C[z_1,z_2]\T S^2\g/V_{2\ta}.
\end{gather*}
Then the conditions (\ref{van}) and (\ref{sym}) are equivalent to
\begin{gather*}
r_0\in z_1z_2(z_1-z_2)^2 P^{\rm sym}\T V_{2\ta},\\
r_1\in z_1z_2(z_1-z_2)P^{\rm sym}\T \Lambda^2\g,\\
r_2\in z_1z_2P^{\rm sym}\T S^2\g/V_{2\ta}.
\end{gather*}
Therefore, the dual quotient space $(V_2/P_+^{\rm sym}V_2)^*$ is isomorphic to $\g**\g$ as a vector space.
It is straightforward to check that the $\g\T\C[t]$ modules structures are also isomorphic.
\end{proof}

We will also need a nonsymmetric version of the  space $V_m/P^{\rm sym}_{m+}V_m$. Namely, let $W_m$ be
a subspace of $z_1\cdots z_m\C[z_1,\dots,z_m]\T \g^{\T m}$ satisfying
\begin{gather}\label{=}
F_{z_i=z_j}\hk \C[z_1,\dots,z_m]|_{z_1=z_2} \T \mu_{i,j} \big(S^2\g/V_{2\ta}\T \g^{\T m-2}\big),\\
\label{'=}
\partial_{z_i}F_{z_i=z_j}\hk \C[z_1,\dots,z_m]|_{z_1=z_2} \T \mu_{i,j} \big(S^2(\g)/V_{2\ta}\T \g^{\T m-2}\big)
\end{gather}
for all $1\le i<j\le m$, where $\mu_{i,j}$ is an operator on the tensor power
$\g^{\T m}$ which inserts the f\/irst two factors on the $i$-th and $j$-th places
respectively:
\begin{gather*}
\mu_{i,j}(v_1\T v_2\T v_2\T \cdots\T v_n)=
v_3\T\cdots\T v_{i-1}\T v_1\T v_i\T \cdots\T v_{j-2}\T v_2\T v_{j-1}\cdots\T v_n.
\end{gather*}
The natural action of the polynomial ring $P_m=\C[z_1,\dots,z_m]$
commutes with the action of the current algebra, def\/ined as in Lemma~\ref{123}. Therefore,
we obtain a structure of $\g\T\C[t]$ module on the quotient space $W_m/P_{m+}W_m$,
where $P_{m+}$ is the ring of positive degree polynomials.
As proved in \cite{FFJMT}, Lemma $5.8$, the symmetric and nonsymmetric constructions
produce the same module (see Proposition~\ref{nsym} below for the precise statement).
We illustrate this in the case $n=2$.

\begin{lemma}\label{W2}
$(W_2/P_{2+}W_2)^*\simeq \g**\g$.
\end{lemma}
\begin{proof}
We use the decomposition $r=r_0+r_1+r_2$ as in Lemma~\ref{n=2} for
$r\in W_m$. The  condi\-tions~(\ref{=}) and (\ref{'=}) are equivalent to the following
conditions on $r_i$:
\begin{gather*}
r_0\in (z_1-z_2)^2 z_1z_2\C[z_1,z_2]\T V_{2\ta},\\
r_1\in (z_1-z_2)z_1z_2\C[z_1,z_2]\T \Lambda^2\g,\\
r_2\in z_1z_2\C[z_1,z_2]\T S^2\g/V_{2\ta}.
\end{gather*}
After passing to the quotient with respect to the action of the algebra $P_{2+}$
we arrive at the isomorphism of vector spaces $(W_2/P_{2+}W_2)^*\simeq \g*\g$.
It is straightforward to check that this is
an isomorphism of $\g\T\C[t]$ modules.
\end{proof}

The following Proposition is proved in \cite[Lemma~5.8]{FFJMT} for $\g=\slt$.
\begin{proposition}\label{nsym}
There exists an isomorphism of $\g\T\C[t]$-modules
\[
V_m/P_{m+}^{\rm sym}V_m\simeq W_m/P_{m+}W_m.
\]
\end{proposition}
\begin{proof}
The proof of the general case dif\/fers from the one from \cite{FFJMT}
by the replacement of the decomposition of the tensor square
of the adjoint representation of $\slt$ by the general decomposition~(\ref{decomp}).
\end{proof}

\begin{remark}
We note that the def\/inition of the space $W_m$ is a bit more involved than
the def\/inition of $V_m$. The reason is that the polynomials used
to construct $W_m$ are not symmetric. In particular, this forces to add the condition
(\ref{'=}) in order to get the isomorphism
$V_m/P_{m+}^{\rm sym}V_m\simeq W_m/P_{m+}W_m$.
\end{remark}

\begin{proposition}
The module $W_m/P_+W_m$ is cocyclic with a cocyclic vector being the class of
\[
r_m=z_1\cdots z_m\prod_{1\le i<j\le m} (z_i-z_j)^2\T (e_\ta)^{\T m}.
\]
\end{proposition}
\begin{proof}
We f\/irst show that if $r\in W_m$ is a nonzero element satisfying
$r\in V_{m\ta}\T \C[z_1,\dots,z_m]$ then either $r\in P_{m+}W_m$ or
there exists an element of the universal enveloping
algebra of $\g\T\C[t]$ which sends $r$ to $r_m$ (here we embed $V_{m\ta}$ into
$\g^{\T m}$ as an irreducible component containing~$(e_\ta)^{\T m}$).
In fact, from conditions (\ref{=}), (\ref{'=}) and the assumption
$r\in V_{m\ta}\T \C[z_1,\dots,z_m]$ we obtain that
$r$ is divisible by the product $\prod_{1\le i<j\le m} (z_i-z_j)^2$.
If $r\notin P_{m+}W_m$ then we obtain
\[
r=x\T \prod_{1\le i<j\le m} (z_i-z_j)^2
\]
with some $x\in V_{m\ta}$. Since $V_{m\ta}$ is irreducible the $\U(\g)$ orbit
of $x$ contains $(e_\ta)^{\T m}$.

Thus it suf\/f\/ices to show that any element $r\in W_m/P_+W_m$ is contained in the $\U(\g\T\C[t])$ orbit
of the image of $V_{m\ta}\T\C[z_1,\dots,z_m]$ in the quotient space.
We prove that if $(a\T t)r=0$ for all $a\in\g$ then $r\in V_{m\ta}\T P_m$.
This would imply the previous statement since
$W_m/P_{m+}W_m$ is f\/inite-dimensional.

So let $r\in W_m$ be some element and $\bar r\in W_m/P_{m+}W_m$ be its class.
Assume that $(a\T t)r=0$ for all $a\in\g$.
We f\/irst show that
$r\in V_{2\ta}\T\g^{\T n-2}\T P_m$.
In fact, there exists a polynomial $l\in P_m$ such that the following holds in
$W_m$:
\[
(a\T t) r = l(z_1,\dots,z_m)r_1
\]
for some $r_1\in W_m$. Consider the space of functions $W^{1,2}_m$ which satisfy conditions
(\ref{=}) and (\ref{'=}) only for $i=1$, $j=2$. (In particular, $W_m\hk W_m^{1,2}$).
We have an isomorphism of $\g\T\C[t]$ modules
\[
W^{1,2}_m\simeq W_2\T \g^{\T m-2}\T \C[z_3,\dots,z_m].
\]
The equality $(a\T t) r=lr_1$ gives
\[
\big(a^{(1)}\T z_1+a^{(2)}\T z_2\big)r=lr_1-\left(\sum_{i=3}^m a^{(i)}\T z_i\right)r.
\]
Using Lemma~\ref{W2}, we obtain
$r\in V_{2\ta}\T \g^{\T n-2}\T P_m$.

The same procedure can be done for all pairs $1\le i<j\le n$.
Now our proposition follows from the following equality  in $\g^{\T n}$:
\begin{gather*}
\bigcap_{1\le i<j\le n} \mu_{i,j} V_{2\ta}\T\g^{\T n-2}=V_{n\ta}.\tag*{\qed}
\end{gather*}\renewcommand{\qed}{}
\end{proof}

The dual module $(W_m/P_m W_m)^*$ is cyclic. We denote by $r'_m$ a cyclic vector which
corresponds to the cocyclic vector $r_m$.

\begin{proposition}\label{surj}
There exists a surjective homomorphism of $\g\T\C[t]$ modules
\begin{gather}\label{gm}
\g^{**m}\to (W_m/P_{m+}W_m)^*
\end{gather}
sending a cyclic vector $e_\ta^{\T m}$ to $r'_m$.
\end{proposition}
\begin{proof}
A relation in the fusion product means that some expression of the form
\begin{gather}
\label{expr}
\sum \al_{i_1\dots i_s} \big(x_1\T t^{i_1}\big)\cdots \big(x_s\T t^{i_s}\big)
\end{gather}
with f\/ixed $t$ degree can be expressed in $\g^{\T m}$ via a linear combination of monomials of
lower $t$-degree.
The coef\/f\/icients of the expression of (\ref{expr}) in terms of the lower degree monomials
are polynomials in evaluation parameters $c_1,\dots, c_n$. Therefore, by the very def\/inition of the action of $P_m$
on $W_m$, the operator of the form (\ref{expr}), which is a relation in the fusion product,
vanishes on $(W_m/P_{m+}W_m)^*$.
\end{proof}

Note  that $W_m$ and $V_m$ are naturally graded by the degree grading on $P_m$ ($P_m^{\rm sym}$).
This grading def\/ines a graded character of the space $V_m/P^{\rm sym}_{m+}$.

\begin{corollary}\label{func}
$\ch_q V_m/P_{m+} V_m\le q^{m^2} \ov{\ch}_{q^{-1}} \g^{**m}$.
\end{corollary}
\begin{proof}
Follows from Proposition~\ref{surj} and Proposition~\ref{nsym}. Note that the factor
$m^2$ is a degree of the cyclic vector $r'_m$.
\end{proof}

\begin{remark}
In the next section we combine Corollary \ref{func} with Corollary \ref{Fs} in
order to compute the character of $L^{\rm gr}$.
\end{remark}

\section{Proofs of the main statements}\label{main}
\begin{proposition}\label{le}
$\ch_q (L_m)^* \le \ch_q V_m\le q^{m^2} \ov{\ch}_{q^{-1}} \g^{**m}/(q)_m.$
\end{proposition}
\begin{proof}
We know that
\[
\ch_q L_m^*\le \ch (L'_m)^*=\ch_q V_m.
\]
 Because of the surjection
(\ref{gm}), the character of $V_m^*$ is smaller than or equal to
$q^{m^2}\ov{\ch}_{q^{-1}} \g^{**m}/(q)_m$ (since $1/(q)_m$ is the character of the space of symmetric
polynomials in $m$ variables). This proves the Proposition.
\end{proof}

\begin{theorem}\label{'}
$L_m\simeq L_m'$.
\end{theorem}
\begin{proof}
Corollary~\ref{Fs} provides an inequality
\[
\ch_q L_m\ge \frac{q^{m^2}}{(q)_m} \ov{\ch}_{q^{-1}} \g^{**m}.
\]
Now from Proposition~\ref{le} we obtain
\[
\frac{q^{m^2}}{(q)_m} \ov{\ch}_{q^{-1}} \g^{**m}\le \ch_q L_m\le \ch_q L_m'\le \frac{q^{m^2}}{(q)_m} \ov{\ch}_{q^{-1}} \g^{**m}.
\]
Theorem is proved.
\end{proof}

\begin{corollary}\label{g*m}
The dual module $(V_m/P^{\rm sym}_{m+}V_m)^*$ and $\g^{**m}$ are isomorphic as $\g\T\C[t]$ modules.
\end{corollary}
\begin{proof}
Follows from Propositions~\ref{surj}, \ref{nsym} and Theorem \ref{'}.
\end{proof}

\begin{corollary}\label{rel}
We have an isomorphism of $\g^{ab}\T\C[t^{-1}]$ modules
\[
L^{\rm gr}\simeq \U\big(\g^{ab}\T\C[t^{-1}]\big)/I,
\]
where $I$ is the minimal $\g$ invariant ideal containing the coefficients of the
current $e_{\ta}(z)^2=0$.
\end{corollary}
\begin{remark}
Corollary \ref{rel} is a generalization of the  $\slt$ case from \cite{FFJMT}. It also
proves a level~1 conjecture from~\cite{F2}.
\end{remark}

\begin{corollary}
The action of the polynomial ring $P_m^{\rm sym}$ on $V_m$ is free.
\end{corollary}
\begin{proof}
Follows from the isomorphism $(L_m/P_{m+}^{\rm sym})^*\simeq \g^{**m}$ and the character equality
\begin{gather*}
\ch_q L_m=\frac{q^{m^2}}{(q)_m} \ov{\ch}_{q^{-1}} \g^{**m}=q^{m^2}\ch_q P_m^{\rm sym} \cdot
\ov{\ch}_{q^{-1}}\g^{**m}.\tag*{\qed}
\end{gather*}  \renewcommand{\qed}{}
\end{proof}

Recall the vectors
\[
w_{i_1,\dots,i_m}=\big(e_\ta\T t^{i_1}\big)\cdots \big(e_\ta\T t^{i_m}\big)v_0\in L.
\]
Let $\bar w_{i_1,\dots,i_m}\in L_m$  be the images of these vectors.

\begin{corollary}
$L_m$ is generated by the set of vectors
\begin{gather}\label{wbar}
\bar w_{i_1,\dots,i_m}, \qquad -m\ge i_1\ge \cdots\ge i_m
\end{gather}
with the action of the algebra $\U(\g\T\C[t])$.
\end{corollary}
\begin{proof}
Recall the element
\[
r_m=z_1\cdots z_m\prod_{1\le i<j\le m} (z_i-z_j)^2\T e_\ta^{\T N}\in V_m.
\]
Since $r_m$ is a cocyclic vector of $V_m/P^{\rm sym}_{m+}$ and the polynomial algebra
acts freely on $V_m$, we obtain that the space
$P^{\rm sym}_{m+}r_m$ is a cocyclic subspace of $V_m$, which means that
for any vector $v\in V_m$ the space $\U(\g\T\C[t])\cdot v$ has a nontrivial
intersection with $P^{\rm sym}_{m+}r_m$. We note that
$P^{\rm sym}_{m+}r_m$ coincides with the subspace of $V_m$ of $\g$ weight $m\ta$.
Dualizing this construction we obtain that the subspace of $L_m$ of $\g$ weight
$-m\ta$ is $\U(\g\T\C[t])$ cyclic. This space is linearly spanned  by the set
(\ref{wbar}). Corollary is proved.
\end{proof}

\begin{corollary}\label{mN}
The space $\gr_m F_\bullet(N)$ is generated by the set of vectors
\begin{gather}\label{-m-N}
\bar w_{i_1,\dots,i_m}, \qquad -m\ge i_1\ge \cdots\ge i_m\ge -N
\end{gather}
with the action of the algebra $\U(\g\T\C[t])$.
\end{corollary}
\begin{proof}
Introduce an increasing  f\/iltration $J_\bullet$ on $L_m$:
\[
J_n=\U(\g\T\C[t])\cdot\mathrm{span}\big\{\bar w_{i_1,\dots,i_m}:\ -i_1-\cdots -i_m\le m^2+n \big\}.
\]
Corollary~\ref{g*m} provides an isomorphism of $\g\T\C[t]$ modules
\[
\U(\g\T\C[t])\cdot \bar w_{i_1,\dots,i_m}/\left(J_{-i_1-\cdots-i_m-1}\cap  \U(\g\T\C[t])\cdot \bar w_{i_1,\dots,i_m}\right)\simeq \g^{**m}.
\]
Thus, since all vectors of the form (\ref{-m-N}) belong to $F_m(N)$, we obtain
\[
\dim \gr_m F_\bullet(N)\ge \binom{N}{m}(\dim\g)^m.
\]
Because of the equality
\[
\dim D(N)=(1+\dim\g)^N=\sum_{m=0}^N \binom{N}{m}(\dim\g)^m,
\]
we conclude that
$\dim \gr_mF_\bullet(N)= \binom{N}{m}(\dim\g)^m$
and hence the whole space $\gr_m F_\bullet(N)$ is gene\-ra\-ted by the vectors
(\ref{-m-N}).
\end{proof}

\begin{proposition}
The induced PBW filtration $F_\bullet(N)\hk D(N)$ coincides with the $t^N$-filtration~$G^\bullet$, i.e.\ $F_m(N)=G^{N-m}$.
\end{proposition}
\begin{proof}
Recall (see Lemma~\ref{GF}) that $G^{N-m}\hk F_m(N)$. Since $v_0$ is proportional to
$(f_\ta\T t^N)^N v_N$,
we obtain that
$\bar w_{i_1,\dots,i_m}\in G^{N-m}$ for
$-m\ge i_1\ge \dots\ge i_m\ge -N$. Therefore, Corollary~\ref{mN} gives
$F_m(N)\hk G^{N-m}$. Proposition is proved.
\end{proof}

\begin{corollary}
The graded component $\gr_m F_\bullet(N)$
is filtered by  $\binom{N}{m}$ copies of $\g^{**m}$.
The character of the space of cyclic vectors of those fusions is equal
to $q^{m^2}\qbin{N}{m}$.
\end{corollary}

We summarize all above in the following theorem:
\begin{theorem}
Let $F_\bullet$ be the PBW filtration on the level one vacuum $\gh$ module $L$.
Then
\begin{enumerate}\itemsep=0pt
\item[\rm a)] $\gr_m F_\bullet$ is filtered by the fusion modules $\g^{**m}$.
\item[\rm b)] The character of the space of cyclic vectors of these $\g^{**m}$ is equal to
$\frac{q^{m^2}}{(q)_m}$.
\item[\rm c)] The induced PBW filtration on Demazure modules $D(N)$ coincides with the double
fusion filtration coming from Theorem~{\rm \ref{FL}}.
\item[\rm d)] The defining relation in $L^{\rm gr}$ is $e_\ta(z)^2=0$.
\end{enumerate}
\end{theorem}

\pdfbookmark[1]{A list of the main notations}{notation}

\section*{A list of the main notations}
\begin{enumerate}\itemsep=0pt
\item[] $\g$ -- simple f\/inite-dimensional algebra;
\item[] $\gh$ -- corresponding af\/f\/ine Kac--Moody algebra;
\item[] $\theta$ -- highest weight  of the adjoint representation of $\g$;
\item[] $e_\ta, f_\ta\in\g$ -- highest and lowest weight vectors of the adjoint representation;
\item[] $L$ -- the basic (vacuum level one) representation of $\gh$;
\item[] $v_0\in L$ -- a highest weight vector of $L$;
\item[] $v_N\in L$ -- an extremal vector of the weight $N\ta$;
\item[] $D(N)\hk L$ -- Demazure module with a cyclic vector $v_N$;
\item[] $\gr_m A_\bullet$ -- $m$-th graded component of the associated graded space with respect
to the f\/iltration $A_\bullet$;
\item[] $F_\bullet$ --  (increasing) PBW f\/iltration on $L$;
\item[] $F_\bullet(N)=F_\bullet\cap D(N)$ -- induced PBW f\/iltration on $D(N)$;
\item[] $G^\bullet$ -- (decreasing) $t^N$-f\/iltration on $D(N)$;
\item[] $V_1*\dots *V_N$ -- associated graded $\g\T\C[t]$ module with respect to the fusion f\/iltration
on the tensor product of cyclic $\g\T\C[t]$ modules;
\item[] $V_1**\dots **V_N$ -- associated graded $\g\T\C[t]$ module with respect to the fusion f\/iltration
on the tensor product of cyclic $\g$ modules;
\item[] $\gr_m (V_1 *\dots *V_N)$ -- $m$-th graded component with respect to the fusion f\/iltration;
\item[] $\ch_q$ -- a graded character def\/ined by the operator $d$;
\item[] $\ov{\ch}_q$ -- a graded character def\/ined by the operator $\bar d$.
\end{enumerate}

\subsection*{Acknowledgements}
EF thanks B.~Feigin and P.~Littelmann for useful discussions.
This work was partially supported by the RFBR  Grants
06-01-00037, 07-02-00799 and NSh-3472.2008.2, by Pierre Deligne fund based on his 2004
Balzan prize  in mathematics, by Euler foundation and by Alexander von Humboldt
Fellowship.


\pdfbookmark[1]{References}{ref}
\LastPageEnding

\end{document}